\documentclass[12pt]{amsart}
\usepackage{amsfonts, amssymb, amsmath, amsthm}
\usepackage{amscd}
\usepackage{color}
\usepackage{enumerate}

\newcommand{\Z}{\mathbb{Z}}

\newcommand{\N}{\mathbb{N}}
\newcommand{\Q}{\mathbb{Q}}
\newcommand{\R}{\mathbb{R}}

\newcommand{\Rcal}{\mathcal{R}}

\newtheorem{theorem}{Theorem}[section]
                            {\end{enumerate}}
\input xy
\xyoption{all}
\numberwithin{equation}{section}
\begin{document}
\keywords{algebraic K-theory, Nil group, finite group ring, strong shift equivalence}
\subjclass[2010]{Primary 19M05; Secondary 37B10}

\title[K Theory Examples]{Explicit examples in $NK_{1}$}

\author{Scott Schmieding}
\address{Department of Mathematics, University of Maryland, College
  Park MD 20742 U.S.A.}
\email{schmiedi@math.umd.edu}

\begin{abstract}
For certain rings $\Rcal$,
we construct explicit matrices representing
nonzero classes in the
algebraic $K$ theory group $NK_{1}(\Rcal)$.
\end{abstract}

\maketitle
%\tableofcontents
\section{Introduction}
The purpose of this note is simply to
exhibit explicit matrices representing non-zero classes in
the algebraic $K$-theory group $NK_{1}(\Rcal)$ (and
thereby in $Nil_{0}(\Rcal) $), for
some rings $\Rcal$ for which $NK_1(\Rcal )$ arises as an obstruction in
\cite{BoSc2} and \cite{BoSc3}. In   \cite{BoSc2}, $\Rcal$ is the integral group ring
of a finite group $G$; our example is for
$G= \Z/4\Z$. (Proof of the nontriviality of $NK_1(\Z [\Z/4])$ can be found in \cite{Weibel2009} or
\cite{SchmiedingNK1FinAbGrps}.) In \cite{BoSc3}, $\Rcal$ is a subring of $\R$;
our example is for $\Rcal = \Q [t^2,t^3,z,z^{-1}]$
(which has many embeddings into $\R$).\\
\indent It seems to be difficult to locate explicit examples
of this sort in the literature.  The arguments to
follow give some indication as to why that might be the case.
The arguments are elementary, and only require carefully tracing through standard arguments and constructions in algebraic K-theory.
However, the actual computation  becomes lengthy,
and leads to fairly large matrix examples.

The computation might be of interest for someone
new to $K$-theory, as an example of the complication buried in
certain exact sequence arguments.\\
\indent The author would like to thank Mike Boyle for helpful comments and suggestions.

\section{Setup}

\indent We will require little setup, all of which can be found in \cite{WeibelBook}, or \cite{Rosenberg1994}.  Always,  $\Rcal$ is an associative ring with $1$.
We consider $K_{1}(\Rcal)$ as the group $GL(\Rcal)/El(\Rcal)$, where $GL(\Rcal)$ denotes the stabilized general linear group over $\Rcal$, and $El(\Rcal)$ the stabilized elementary subgroup of $GL(\Rcal)$. There is a map $K_{1}(\Rcal[t]) \to K_{1}(\Rcal)$ induced by $t \mapsto 0$, and the kernel of this map is defined to be $NK_{1}(\Rcal)$.
Higman's trick implies $NK_{1}(\Rcal)$
is the set of elements of $K_{1}(\Rcal[t])$ which contain a matrix of the form
$I-tN$, with $N$ a nilpotent matrix over $\Rcal$.

The group
$Nil_{0}(\Rcal)$ is defined from
\textbf{Nil}$\Rcal$, the nilpotent category over $\Rcal$.
The objects of this category
are  pairs $(P,f)$, where $P$ is a finitely generated projective $\Rcal$-module
and $f$ is a nilpotent endomorphism of $P$.
A morphism $(P,f) \to (Q,g)$ is an $\Rcal$ module homomorphism
$h: P\to Q$ such that $hf=gh$.
\textbf{Nil}$\Rcal$ acquires an exact structure via
the forgetful functor \textbf{Nil}$\Rcal \to $ \textbf{Proj}$\Rcal$ given by $(P,f) \mapsto P$:  a sequence in \textbf{Nil}$\Rcal$ is exact if its image under this forgetful functor is exact. One may then consider the $K$-group $K_{0}($\textbf{Nil}$\Rcal)$ of the exact category $\textbf{Nil}\Rcal$
(see \cite[II.7]{WeibelBook}).
The cokernel of the map
$K_{0}(\textbf{Proj}\Rcal) \to K_{0}(\textbf{Nil}\Rcal) $
given by $[P]\mapsto [(P,0)]$
is denoted $Nil_{0}(\Rcal)$.
A well-known isomorphism $NK_{1}(\Rcal ) \to
Nil_{0}(\Rcal)$ is induced by $N\mapsto I-tN$, where $N$ denotes a nilpotent
matrix over $\Rcal$  (viewed as an endomorphism of $\mathcal R^n$, where
$N$ is $n\times n$).
\\

\indent For $k \ge 1$, we recall an important endomorphism of $NK_{1}(\Rcal)$. The Verschiebung map $V_{k} $
acts on $NK_{1}(\Rcal)$ via $V_{k} ([1-tN]) = [1-t^{k}N]$, and acts on $Nil_{0}(\Rcal)$ via
$$V_{k}([N]) =
\Bigg[
\left(\begin{smallmatrix} 0 & 0 & \cdots & 0 & N \\ 1 & 0 & 0 & \cdots & 0 \\ 0 & 1 & 0 & \cdots & 0 \\ \vdots & 0 & \ddots & \ddots & \vdots \\ 0 & \cdots & 0 & 1 & 0 \end{smallmatrix} \right)
\Bigg] \ . $$
For $k \ge 1$ there are associated endomorphisms called the Frobenius maps $F_{k}$, which act on $NK_{1}(\Rcal)$ via
$F_{k} ([I-tN]) = [I-tN^{k}]$,  and
acts on $Nil_{0}(\Rcal)$ via $F_{k}([N]) = [N^{k}]$. We will only require the Verschiebung map in the present note.\\

\section{A matrix nontrivial in $NK_{1}(\mathbb{Q}[t^{2},t^{3},z,z^{-1}])$}
Consider the ring $\mathbb{Q}[t^{2},t^{3},z,z^{-1}]$. For consistency, we let $NK_{1}(\mathbb{Q}[t^{2},t^{3},z,z^{-1}])$ denote the kernel of the map $K_{1}(\mathbb{Q}[t^{2},t^{3},z,z^{-1},s]) \stackrel{s \mapsto 0}\to K_{1}(\mathbb{Q}[t^{2},t^{3},z,z^{-1}])$. In this section, we show the following.
\begin{theorem}
\label{Example1}
\begin{enumerate}
\item
The class of the matrix
$$\begin{pmatrix} 1-(1+z^{-1})s^{4}t^{4} & (z-1)(s^{2}t^{2}-s^{3}t^{3}) \\ (1-z^{-1})(s^{2}t^{2})(1+st+s^{2}t^{2}+s^{3}t^{3}) & 1+(z-1)(s^{4}t^{4}) \end{pmatrix}$$
is not zero in $NK_{1}(\mathbb{Q}[t^{2},t^{3},z,z^{-1}])$. \\
\item
\indent The class of the matrix
$$\begin{pmatrix}
0&0 &         0&(1-z)t^2 &     0&(1-z)(-t^3) & (1-z^{-1})t^4&0 & 0&0\\
0&0 & (z^{-1}-1)t^2&0 & (z^{-1}-1)t^3&0 & (z^{-1}-1)t^4&(1-z)t^4 & (z^{-1}-1)t^5&0\\
1&0 & 0&0 & 0&0 & 0&0 &0&0 \\
0&1 & 0&0 & 0&0 & 0&0 &0&0 \\
0&0 & 1&0 & 0&0 & 0&0 &0&0 \\
0&0 & 0&1 & 0&0 & 0&0 &0&0 \\
0&0 & 0&0 & 1&0 & 0&0 &0&0 \\
0&0 & 0&0 & 0&1 & 0&0 &0&0 \\
0&0 & 0&0 & 0&0 & 1&0 &0&0 \\
0&0 & 0&0 & 0&0 & 0&1 &0&0
\end{pmatrix}
$$
is non-zero in $Nil_{0}(\mathbb{Q}[t^{2},t^{3},z,z^{-1}]$.
\end{enumerate}
\end{theorem}

It follows from \cite[Theorem 5.4]{BoSc1} that the nilpotent matrix in $(2)$ of \ref{Example1} is not strong shift equivalent over the ring to the zero matrix (see Section 5). Since this ring may be embedded in $\mathbb{R}$, this provides an example of a matrix over a subring of the reals which is shift equivalent, but not strong shift equivalent, to zero.\\
%\cite[Theorem 5.4]{BoSc1}.
\indent The remainder of this section is devoted to proving
\ref{Example1}. An outline for the construction is as follows: we
begin with a non-zero class in $K_{1}(\mathbb{Q}[t,s]/I)$ which lies
in the kernel of the map induced by $s \mapsto 0$. Such classes are
easy to find. We then proceed by applying a collection of maps to this
element, taking care that at each stage of the composition, the
element remains non-zero, and ending in
$K_{1}(\mathbb{Q}[t^{2},t^{3},z,z^{-1},s])$. The resulting element
still lies in the kernel upon sending $s \mapsto 0$, so the final
element lies in $NK_{1}[\mathbb{Q}[t^{2},t^{3},z,z^{-1}])$. The maps
which will be applied are shown in the diagram below, starting bottom left and ending top right:
%\annotation{See changes at the vertical arrow (color didn't work
%  there).}
\[ \xymatrix{
 & K_{0}(\mathbb{Q}[t^{2},t^{3},s],I)
%\ar[d]
\ar[r]^{(p_{2})_{*}} & K_{0}(\mathbb{Q}[t^{2},t^{3},s])
\ar[r]^{\hspace{-.2in}\cdot z} &
K_{1}(\mathbb{Q}[t^{2},t^{3},z,z^{-1},s]) \\ K_{1}(\mathbb{Q}[t,s]/I)
\ar[r]^{\partial} &
 K_{0}(\mathbb{Q}[t,s],I)
\ar[u]^{\phi_{*}}_{\cong} & & }  \\
\]
\\
\indent Consider the ideal $I = t^{2}\mathbb{Q}[t] \subset \mathbb{Q}[t,s]$. There is an exact sequence
$$0 \to I \to \mathbb{Q}[t,s] \to \mathbb{Q}[t,s]/I \to 0$$
which yields an exact sequence in $K$-theory (see \cite[2.5.4]{Rosenberg1994})
$$\cdots \to  K_{1}(\mathbb{Q}[t,s]) \stackrel{^{\pi_{*}}}\to K_{1}(\mathbb{Q}[t,s] / I) \stackrel{^{\partial}}\to K_{0}(\mathbb{Q}[t,s],I) \stackrel{(p_{2})_{*}}\to K_{0}(\mathbb{Q}[t,s]) \stackrel{\pi_{*}}\to K_{0}(\mathbb{Q}[t,s] / I)$$
where $K_{i}(\mathbb{Q}[t,s],I)$ refers to the relative groups defined via the double
$$D(\mathbb{Q}[t,s],I) = \{(x,y) \in \mathbb{Q}[t,s] \times \mathbb{Q}[t,s] \hspace{.01in} | \hspace{.01in} x-y \in I\}$$ of $\mathbb{Q}[t,s]$ along $I$ (see \cite[1.5.3]{Rosenberg1994}), $\pi_{*}$ is induced by $\pi:\mathbb{Q}[t,s] \to \mathbb{Q}[t,s] / I$, $(p_{2})_{*}$ is induced by the projection onto the second coordinate $p_{2}:D(\mathbb{Q}[t,s],I) \to \mathbb{Q}[t,s]$, and $\partial$ is the boundary map. Consider the class $1+ts \in K_{1}(\mathbb{Q}[t,s] / I)$. \\

\indent \textbf{Step 1 - computing $\partial$:}  Since $1+ts$ is not in the image of $\pi_{*}:K_{1}(\mathbb{Q}[t,s]) \to K_{1}(\mathbb{Q}[t,s] / I)$, $\partial(1+ts)$ represents a non-zero class in $K_{0}(\mathbb{Q}[t,s],I)$. We proceed by computing $\partial(1+ts)$ explicitly, which is done via a standard clutching construction. An outline of this can be found in \cite[2.5.4]{Rosenberg1994}.\\
\indent First consider $M_{1+ts} = \{(x,y) \in \mathbb{Q}[t,s]^{2} \hspace{.01in} | \hspace{.01in} y-x(1+ts) \in I\}$ (thinking of $x,y$ as row vectors). This is a projective $D(\mathbb{Q}[t,s],I)$-module, and we get $\partial(1+ts)=[M_{1+ts}] - [D(\mathbb{Q}[t,s],I)]$. For computational purposes it turns out to be more useful to compute the class of idempotent matrices representing $[M_{1+ts}]$ and $[D(\mathbb{Q}[t,s],I)]$. For this, first note the product

\begin{align*}
A&=\begin{pmatrix} 1+st+s^{2}t^{2}+s^{3}t^{3} & -s^{2}t^{2} \\
  s^{2}t^{2} & 1-st \end{pmatrix}\\
& = \begin{pmatrix} 1 & 1+st \\ 0 & 1\end{pmatrix} \begin{pmatrix} 1 &
  0 \\ -(1+st) & 1 \end{pmatrix} \begin{pmatrix}1 & 1+st \\ 0 &
  1 \end{pmatrix} \begin{pmatrix} 0 & -1 \\ 1 & 0 \end{pmatrix}
\end{align*}
%\annotation{I changed the displays only in the display (not in any
%  formula), using align*}
 is a lift of $\begin{pmatrix} 1+st & 0 \\ 0 & 1-st \end{pmatrix}$ (so $\pi(A) = \begin{pmatrix} 1+st & 0 \\ 0 & 1-st \end{pmatrix}$), and there is an isomorphism $$j:M_{1+st} \oplus M_{1-st} \to D(\mathbb{Q}[t,s],I)^{2}$$
given by

\begin{align*} j:\begin{pmatrix} (x,y) & (u,v) \end{pmatrix}
&\to
  (\begin{pmatrix}x & u \end{pmatrix}A,\begin{pmatrix} y &
    v \end{pmatrix}) \\
&\to \begin{pmatrix}(\pi_{1}(\begin{pmatrix}x &
      u \end{pmatrix}A),\pi_{1}(\begin{pmatrix} y & v \end{pmatrix}))
    & (\pi_{2}(\begin{pmatrix}x &
      u \end{pmatrix}A),\pi_{2}(\begin{pmatrix} y &
      v \end{pmatrix})) \end{pmatrix} \\
&= \begin{pmatrix}(\pi_{1}(\begin{pmatrix}x & u \end{pmatrix}A),y) &
  (\pi_{2}(\begin{pmatrix}x & u \end{pmatrix}A),v) \end{pmatrix}
\end{align*}
where $\pi_{1}$ and $\pi_{2}$ denote projection on to the 1st and 2nd component, respectively. Using this isomorphism, one can check that
\[ \xymatrix{
& M_{1+st} \oplus M_{1-st} \ar[r]^{j} \ar[d]^{id \oplus 0} & D(\mathbb{Q}[t,s],I)^{2} \ar[d]^{B}\\
& M_{1+st} \oplus M_{1-st} \ar[r]^{j} & D(\mathbb{Q}[t,s],I)^{2} }  \\
\]
commutes, where $B = \begin{pmatrix} (1-s^{4}t^{4},1) & ((-s^{2}t^{2})(1+st+s^{2}t^{2}+s^{3}t^{3}),0) \\ (s^{3}t^{3}-s^{2}t^{2},0) & (s^{4}t^{4},0) \end{pmatrix}$. Thus, the idempotent matrix $B$ represents the class of $M_{1+ts}$ in $K_{0}(\mathbb{Q}[t,s],I)$. For ease of notation in what follows, we will express $B$ instead as an ordered pair of matrices $B = (B_{1},B_{2})$ coming from each of the components of the entries of $B$, so
$$B = (\begin{pmatrix}1-s^{4}t^{4} & (-s^{2}t^{2})(1+st+s^{2}t^{2}+s^{3}t^{3}) \\ s^{3}t^{3}-s^{2}t^{2} & s^{4}t^{4} \end{pmatrix}, \begin{pmatrix} 1 & 0 \\ 0 & 0 \end{pmatrix} )$$
and
$$B_{1} = \begin{pmatrix}1-s^{4}t^{4} & (-s^{2}t^{2})(1+st+s^{2}t^{2}+s^{3}t^{3}) \\ s^{3}t^{3}-s^{2}t^{2} & s^{4}t^{4} \end{pmatrix}, B_{2} = \begin{pmatrix} 1 & 0 \\ 0 & 0 \end{pmatrix}$$
Since $[D(\mathbb{Q}[t,s],I)] = [\begin{pmatrix}1 & 0 \\ 0 & 0 \end{pmatrix},\begin{pmatrix} 1 & 0 \\ 0 & 0 \end{pmatrix}]$ in idempotent form, we get $\partial(1+ts) = [B] - [\begin{pmatrix}1 & 0 \\ 0 & 0 \end{pmatrix},\begin{pmatrix} 1 & 0 \\ 0 & 0 \end{pmatrix}]$. For notational ease, let $P = \begin{pmatrix} 1 & 0 \\ 0 & 0 \end{pmatrix}$. \\

\indent \textbf{Step 2 - Excision isomorphism $\phi_{*}$:} Since $K_{0}$ has the excision property, there are isomorphisms $K_{0}(\mathbb{Q}[t,s],I) \stackrel{\gamma_{1,*}^{-1}}\cong K_{0}(I) \stackrel{\gamma_{2,*}}\cong K_{0}(\mathbb{Q}[t^{2},t^{3},s],I)$. Let
$$\phi_{*} = \gamma_{2,*} \circ \gamma_{1,*}^{-1} :K_{0}(\mathbb{Q}[t,s],I) \to K_{0}(\mathbb{Q}[t^{2},t^{3},s],I)$$
denote the composition of these isomorphisms. Here $K_{0}(I)$ denotes $K_{0}$ of the non-unital ring $I$ \cite[1.5.7]{Rosenberg1994}. This is defined by unitizing $I$, i.e. forming the ring $I_{+} = I \oplus \mathbb{Z}$ with multiplication given by $(x,n) \cdot (y,m) = (xy + ny + mx,mn)$, and defining $K_{0}(I)$ to be the kernel of the induced map from the surjection on to the second factor, $K_{0}(I) = ker(K_{0}(I_{+}) \to K_{0}(\mathbb{Z}))$. The isomorphism $\gamma_{1,*}:K_{0}(I) \to K_{0}(\mathbb{Q}[t,s],I)$ is induced by the map $(x,n) \to (n,n+x)$. Letting $e_{2} = (A^{T})^{-1}DA^{T} = \begin{pmatrix}1-s^{4}t^{4} & s^{2}t^{2} - s^{3}t^{3} \\ s^{2}t^{2}(1+st+s^{2}t^{3}+s^{3}t^{3}) & s^{4}t^{4} \end{pmatrix}$ (transposes appear because in the isomorphism $j$ above we are acting on row vectors), a computation (see \cite[1.5.9]{Rosenberg1994} for details) gives
\belowdisplayskip=0pt
$$\gamma_{*}^{-1}([B] - [P,P]) = [e_{2}-P,P]-[\begin{pmatrix}0 & 0 \\ 0 & 0 \end{pmatrix},P]$$
Applying $\gamma_{2,*}$ gives
$$\gamma_{2,*}([e_{2}-P,P]-[\begin{pmatrix}0 & 0 \\ 0 & 0 \end{pmatrix},P]) = [P,e_{2}]-[P,P]$$
so
$$\phi_{*}([B]-[P,P]) = [P,e_{2}]-[P,P] \in K_{0}(\mathbb{Q}[t^{2},t^{3},s],I)$$
\\
\textbf{Step 3 - computing $(p_{2})_{*}$:} The map $(p_{2})_{*}:K_{0}(\mathbb{Q}[t^{2},t^{3},s],I) \to K_{0}(\mathbb{Q}[t^{2},t^{3},s])$ is induced by the projection $p_{2}:D(\mathbb{Q}[t^{2},t^{3},s],I) \to \mathbb{Q}[t^{2},t^{3},s]$ onto the second coordinate. Thus $(p_{2})_{*}([P,e_{2}]-[P,P]) = [e_{2}]-[P]$. We claim this class is non-zero in $K_{0}(\mathbb{Q}[t^{2},t^{3},s])$. To see this, note there is a splitting map $\varphi:\mathbb{Q}[s] \to \mathbb{Q}[t^{2},t^{3},s]$ for $q:\mathbb{Q}[t^{2},t^{3},s] \to \mathbb{Q}[t^{2},t^{3},s]/I = \mathbb{Q}[s]$, which implies that the map $q_{*}:K_{1}(\mathbb{Q}[t^{2},t^{3},s]) \to K_{1}(\mathbb{Q}[t^{2},t^{3},s] / I)$ is surjective. This in turn implies that the boundary map $\partial$ in the exact sequence
$$\cdots \to  K_{1}(\mathbb{Q}[t^{2},t^{3},s]) \stackrel{q_{*}}\to K_{1}(\mathbb{Q}[t^{2},t^{3},s] / I) \stackrel{\partial}\to K_{0}(\mathbb{Q}[t^{2},t^{3},s],I) \stackrel{(p_{2})_{*}}\to K_{0}(\mathbb{Q}[t^{2},t^{3},s]) \to \cdots $$
must be zero. Thus $(p_{2})_{*}$ must be injective. Altogether we have the non-zero class $[e_{2}]-[P] \in K_{0}(\mathbb{Q}[t^{2},t^{3},s],I)$. \\

\indent \textbf{Step 4 - computing $\cdot$ z:} Finally, for any ring
$T$, there is an injective map (see \cite[III.3.5.2]{WeibelBook})
$\cdot z:K_{0}(T) \to K_{1}(T[z,z^{-1}])$ given by $\cdot z:[Q] \to
[I+(z-1)Q]$, where $Q$ is an idempotent matrix over $T$. Thus we apply
this map to the idempotent $e_{2}$ to get $[I+(z-1)e_{2}] \in
K_{1}(\mathbb{Q}[t^{2},t^{3},s,z,z^{-1}])$, and to $P
= \begin{pmatrix} 1 & 0 \\ 0 & 0 \end{pmatrix}$ to get the difference
\begin{align*}
&\ \cdot z ([e_{2}] - [\begin{pmatrix} 1 & 0 \\ 0 & 0 \end{pmatrix}])
\\
 =&\  [I+(z-1)e_{2}] - [\begin{pmatrix} z & 0 \\ 0 & 1 \end{pmatrix}] \\
=&\  [\begin{pmatrix} 1+(z-1)(1-s^{4}t^{4}) &
  (z-1)(s^{2}t^{2}-s^{3}t^{3}) \\
  (z-1)(s^{2}t^{2})(1+st+s^{2}t^{2}+s^{3}t^{3}) &
  1+(z-1)(s^{4}t^{4}) \end{pmatrix}] - [\begin{pmatrix} z & 0 \\ 0 &
  1 \end{pmatrix}]\\
 =&\  [\begin{pmatrix} z^{-1}+(1-z^{-1})(1-s^{4}t^{4}) &
   (z-1)(s^{2}t^{2}-s^{3}t^{3}) \\
   (1-z^{-1})(s^{2}t^{2})(1+st+s^{2}t^{2}+s^{3}t^{3}) &
   1+(z-1)(s^{4}t^{4}) \end{pmatrix}] \\
 =&\  [\begin{pmatrix} 1-(1+z^{-1})s^{4}t^{4} &
   (z-1)(s^{2}t^{2}-s^{3}t^{3}) \\
   (1-z^{-1})(s^{2}t^{2})(1+st+s^{2}t^{2}+s^{3}t^{3}) &
   1+(z-1)(s^{4}t^{4}) \end{pmatrix}]
\end{align*}
in $K_{1}(\mathbb{Q}[t^{2},t^{3},s,z,z^{-1}])$. One can check easily that the above class maps to $[I]$ under the map $s \to 0$, and hence lies in $NK_{1}(\mathbb{Q}[t^{2},t^{3},z,z^{-1},s]$. Lastly, non-triviality of the class was justified at each stage. \\
\indent To find the corresponding class in $Nil_{0}$, let $I-M$ denote this matrix found above, so we have $M$ as
\[
M=\begin{pmatrix}
(1-z^{-1})s^{4}t^{4} & (1-z)(s^{2}t^{2}-s^{3}t^{3})  \\
(z^{-1}-1)(s^{2}t^{2})(1+st+s^{2}t^{2}+s^{3}t^{3}) &
(1-z)(s^{4}t^{4}) \end{pmatrix}
= \sum_{i=1}^5 s^iM_i
\]
with the
$M_i$ over $\mathbb{Q}[t^{2},t^{3},z,z^{-1},s]$. Under the isomorphism $NK_{1} \to Nil_{0}$ we obtain (see \cite{BoSc1}) a nilpotent matrix $N$ over $\mathbb{Q}[t^{2},t^{3},z,z^{-1},s]$,
\begin{align*}
&N\ =\
\begin{pmatrix}
M_1 & M_2 & M_3 & M_4 & M_5 \\
I       &  0    &   0    &  0    &  0    \\
0      &  I     &   0    &  0    &  0    \\
0      &  0    &   I    &  0    &  0    \\
0      &  0    &   0    &  I    &  0
\end{pmatrix} \ =\ \\
&
\begin{pmatrix}
0&0 &         0&(1-z)t^2 &     0&(1-z)(-t^3) & (1-z^{-1})t^4&0 & 0&0\\
0&0 & (z^{-1}-1)t^2&0 & (z^{-1}-1)t^3&0 & (z^{-1}-1)t^4&(1-z)t^4 & (z^{-1}-1)t^5&0\\
1&0 & 0&0 & 0&0 & 0&0 &0&0 \\
0&1 & 0&0 & 0&0 & 0&0 &0&0 \\
0&0 & 1&0 & 0&0 & 0&0 &0&0 \\
0&0 & 0&1 & 0&0 & 0&0 &0&0 \\
0&0 & 0&0 & 1&0 & 0&0 &0&0 \\
0&0 & 0&0 & 0&1 & 0&0 &0&0 \\
0&0 & 0&0 & 0&0 & 1&0 &0&0 \\
0&0 & 0&0 & 0&0 & 0&1 &0&0
\end{pmatrix}
\end{align*}
which is nontrivial as an element of $Nil_{0}(\mathbb{Q}[t^{2},t^{3},z,z^{-1},s])$.\\

\indent One can of course use the above technique to generate many more explicit non-zero classes: simply start with any unit of the form $a + bst$ in $\mathbb{Q}[t,s]/t^{2}\mathbb{Q}[t,s]$, $a,b \in \mathbb{Q}$, and apply the sequence of maps $\cdot z (p_{2})_{*} \phi_{*} \partial$ to $a + bst$ as above.

\section{A matrix nontrivial in $\mathbb{Z}G$, for
$G=\mathbb{Z}/4$}
\indent This section is concerned with the integral group rings of
finite cyclic groups. Let $p$ be a prime, and $\mathbb{Z}/p^{n}$
denote a cyclic group of order $p^{n}$.
In \cite[Theorem 3.12]{MartinThesis} (except for a few cases), and later in \cite{SchmiedingNK1FinAbGrps}, it is shown that
$NK_{1}(\mathbb{Z}[\mathbb{Z}/p^{n}]) \ne 0$ for $n \ge 2$. In \cite{Weibel2009} it is also shown that $NK_{1}(\mathbb{Z}[G])$ is not zero for $G = \mathbb{Z}/4$, along with $G = D_{4}$, the dihedral group. The technique in \cite{SchmiedingNK1FinAbGrps} is an extension of that found in \cite[Theorem1.4]{Weibel2009}, and uses the Milnor square

\begin{center}
 \[ \xymatrix{
 \mathbb{Z}[\mathbb{Z}/p^{n}] \ar[r] \ar[d] & \mathbb{Z}[\zeta_{p^{n}}] \ar[d] \\
 \mathbb{Z}[\mathbb{Z}/p] \ar[r] & \mathbb{Z}[\zeta_{p^{n}}]/(1-\zeta_{p^{n}}^{p}) \\}
\]
\end{center}
where $\zeta_{n}$ denote a primitive $n$th root of unity, and $\mathbb{Z}[\zeta_{n}]$ the ring of integers of $\mathbb{Q}[\zeta_{n}]$. The bottom right term is isomorphic to $\mathbb{Z}_{p}[t]/(t^{p})$ (\cite{SchmiedingNK1FinAbGrps}), and the square yields the Mayer-Vietoris sequence
\begin{equation}\label{diagram}
\begin{gathered}
\xymatrix@C=.8em @R=.7em{
&&NK_{2}(\mathbb{Z}[\mathbb{Z}/p^{n}]) \ar[r] & NK_{2}(\mathbb{Z}[\zeta_{p^{n}}]) \oplus NK_{2}(\mathbb{Z}[\mathbb{Z}/p]) \ar[r] & NK_{2}(\mathbb{Z}_{p}[t]/(t^{p})) \\
%\ar@{->} `r/8pt[d] `/10pt[l] `^dl[ll]|{\alpha^0} `^r/3pt[dll] [dll]
&\ar[r] & NK_{1}(\mathbb{Z}[\mathbb{Z}/p^{n}]) \ar[r] & NK_{1}(\mathbb{Z}[\zeta_{p^{n}}]) \oplus NK_{1}(\mathbb{Z}[\mathbb{Z}/p]) \ar[r] & NK_{1}(\mathbb{Z}_{p}[t]/(t^{p})) \ar[r] & \cdots }\\
\end{gathered}
\end{equation}
Since $\mathbb{Z}[\zeta_{p^{n}}]$ is regular, $NK_{i}(\mathbb{Z}[\zeta_{p^{n}}])=0$ for $i=1,2$, and $NK_{1}(\mathbb{Z}[\mathbb{Z}/p])=0$ from \cite{Harmon1987}. Thus $NK_{1}(\mathbb{Z}[\mathbb{Z}/p^{n}])$ is isomorphic to the cokernel of $NK_{2}(\mathbb{Z}[\mathbb{Z}/p]) \to NK_{2}(\mathbb{Z}_{p}[t]/(t^{p}))$. A presentation of this cokernel is given in \cite{SchmiedingNK1FinAbGrps} using the computation of $NK_{2}(\mathbb{Z}_{p}[t]/(t^{p}))$ by van der Kallen and Stienstra found in \cite{vdKStienstra1984}. \\
\indent We use this method to produce concrete non-zero classes in
$NK_{1}(\mathbb{Z}[\mathbb{Z}/4])$. For the case at hand, namely
$\mathbb{Z}[\mathbb{Z}/4]$, a $p=2$ case of the argument also appears in Weibel \cite[Theorem 1.4]{Weibel2009}, and we mimic the notation found there.\\

\begin{theorem}
\label{Example2}
The class of the matrix $\begin{pmatrix} A & B \\ C & D \end{pmatrix}$ in $NK_{1}(\mathbb{Z}[\mathbb{Z}/4])$
with $$A=1 - (1-\sigma^{2})(x-2x^2+2x^3-\sigma + x\sigma +x^{2} \sigma)$$
$$B=(\sigma^{2}-1)(1+2x-x^{2}-x^{3}-2x^{4}+\sigma-x\sigma -2x^{2}\sigma-3x^{3}\sigma+2x^{4}\sigma) $$
$$C=(\sigma^{2}-1)(-1+2x-5x^{2}+7x^{3}-3x^{4}+2x^{5}-\sigma +2x\sigma -2x^{3}\sigma +3x^{4}\sigma -2x^{5}\sigma) $$
$$D=1 - (1-\sigma^{2})(2+x-2x^{2}-4x^{4}-2x^{5}+\sigma -3x\sigma -x^{2}\sigma -4x^{3}\sigma +6x^{4}\sigma -4x^{5}\sigma +4x^{6}\sigma)$$
is non-zero.
\end{theorem}
To construct a corresponding non-zero class in $Nil_{0}(\mathbb{Z}[\mathbb{Z}/4])$, one could now apply Higman's trick. Since the matrix contains powers of $x$ up to and including 6, this would yield a $12 \times 12$ nilpotent matrix. \\

The remainder of the section is devoted to verifying \ref{Example2}. \\

\indent Let $G = [\mathbb{Z}/4]$ be the cyclic group of order $4$ generated by $\sigma$, and let $\Rcal = \mathbb{Z}[G]$. As described above we have the square \\
\[ \xymatrix{
& & \Rcal \ar[r]^{\sigma \to i} \ar[d]_{\sigma^{2} \to 1} & \mathbb{Z}[i] \ar[d]^{i \to 1+\epsilon} \\
& & \mathbb{Z}[\mathbb{Z}/2] \ar[r]^{q} & \mathbb{F}_{2}[\epsilon]/(\epsilon ^{2})  }  \\
\]
with $q$ the map induced by the quotient $\mathbb{Z} \to \mathbb{F}_{2}$, yielding the Mayer-Vietoris sequence (see \cite{Weibel2009}), a part of which reads\\
$$ \cdots \to NK_{2}(\mathbb{Z}[\mathbb{Z}/2]) \stackrel{q}\to NK_{2}(\mathbb{F}_{2}[\epsilon]/(\epsilon ^{2})) \stackrel{\partial}\to NK_{1}(\Rcal) \to \cdots $$
Furthermore, Lemma 1.2 in \cite{Weibel2009} implies $Im(q) = ker(D: NK_{2}(\mathbb{F}_{2}[\epsilon] /(\epsilon ^{2})) \to \Omega_{\mathbb{F}_{2}[x]}$, where $\Omega_{\mathbb{F}_{2}[x]}$ denotes the K\"{a}hler differentials of $\mathbb{F}_{2}[x]$, and $D$ is the map $D(\langle f \epsilon, g+g^{\prime} \epsilon \rangle) = f \hspace{.03in}dg$. Here $\langle \hspace{.03in}, \rangle$ denotes the Dennis-Stein symbol in $K_{2}$ (see \cite[III.5.11]{WeibelBook}). Thus, choosing for example $\langle \epsilon, x+\epsilon \rangle$, we have $D(\langle \epsilon, x+\epsilon \rangle) = dx \ne 0$, so $\langle \epsilon, x+\epsilon \rangle \not \in Im(q)$. It follows that $\partial(\langle \epsilon, x+\epsilon \rangle) \ne 0$ in $NK_{1}(\Rcal)$. \\
\indent It remains to compute the boundary map $\partial(\langle \epsilon, x+\epsilon \rangle)$. This is obtained (see \cite[III.5.8]{WeibelBook}) from the composition (bottom left to top right)\\
\[ \xymatrix{
& & K_{1}(\Rcal[x],(1-\sigma^{2})) \ar[d]^{\psi}_{\cong} \ar[r]^{\hspace{.3in}j} & K_{1}(\Rcal[x]) \\
& K_{2}(\mathbb{F}_{2}[\epsilon,x]/(\epsilon^{2})) \ar[r]^{\partial_{1}} & K_{1}(\mathbb{Z}[i][x],(2)) & } \\
\]

Here $\psi$ is the induced map from $\Rcal[x] \to \mathbb{Z}[i][x]$,
which is an isomorphism, since the map $\sigma \to i$ takes the ideal
$(1-\sigma^{2})$ isomorphically onto the ideal $(2)$. The map
$\partial_{1}$ is the standard boundary map for the long exact from an
ideal, whose computation is routine\footnote{We include the calculation in an appendix, for completeness.}, yielding\\
$$\partial_{1}(\langle \epsilon, x+\epsilon \rangle) = $$
$$[YZ]$$
where
$$Y:=e_{21}(-x+1-i+(1-i)x^{2})e_{12}(1-i)e_{21}(x+i-1)e_{12}(i-1)$$
$$Z:=e_{12}(1)e_{21}(-1)e_{12}(1)e_{12}((i-1)x-1)e_{21}(1+(i-1)x)e_{12}((i-1)x-1)$$\\
Lifting this up to $K_{1}(\Rcal[x],(1-\sigma^{2}))$ via the vertical isomorphism $\psi$ gives
$$ \begin{pmatrix} A & B \\
C & D \end{pmatrix}$$
with $$A=1 - (1-\sigma^{2})(x-2x^2+2x^3-\sigma + x\sigma +x^{2} \sigma)$$
$$B=(\sigma^{2}-1)(1+2x-x^{2}-x^{3}-2x^{4}+\sigma-x\sigma -2x^{2}\sigma-3x^{3}\sigma+2x^{4}\sigma) $$
$$C=(\sigma^{2}-1)(-1+2x-5x^{2}+7x^{3}-3x^{4}+2x^{5}-\sigma +2x\sigma -2x^{3}\sigma +3x^{4}\sigma -2x^{5}\sigma) $$
$$D=1 - (1-\sigma^{2})(2+x-2x^{2}-4x^{4}-2x^{5}+\sigma -3x\sigma -x^{2}\sigma -4x^{3}\sigma +6x^{4}\sigma -4x^{5}\sigma +4x^{6}\sigma)$$
Applying $j$ yields the class of $\begin{pmatrix} A & B \\ C & D \end{pmatrix}$ in $K_{1}(\Rcal[x])$. This class is our desired element in $NK_{1}(\Rcal)$.\\

\indent The technique above can be repeated, using other symbols $\langle \hspace{.1in}, \hspace{.1in} \rangle$ in $K_{2}(\mathbb{F}_{2}[\epsilon])$, although the computations become quite lengthy.
\section{Strong shift equivalence}

Let $A,B$ be square matrices over a ring $\Rcal$
(not necessarily of the same size).

$A$ and $B$
are {\it elementary strong shift equivalent over $\mathcal R$} (ESSE-$\Rcal$)
if there exist matrices $U,V$ over $\Rcal$ such that $A=UV$ and $B=VU$.
$A$ and $B$
are {\it strong shift equivalent over $\mathcal R$ (SSE-$\Rcal$)} if there
are matrices $A=A_0, A_1, \dots , A_{\ell}=B$ such that for $1\leq i\leq \ell$,
$A_i$ and $A_{i-1}$ are ESSE-$\Rcal$.
$A$ and $B$ are {\it shift equivalent}
over $\Rcal$ if there exist matrices $U,V$ over $\Rcal $ and $\ell $ in $\N$  such that
$A^{\ell} = UV, B^{\ell}=VU, AU = UB$ and $VA=BV$.

It is proved in \cite{BoSc1} that for any nilpotent matrix $N$ over $\Rcal$,
the matrix $A \oplus N$ is SSE-$\Rcal$ to $A$ if and only if $N$ is trivial
as an element of $Nil_0(\Rcal )$.
Moreover, if $B$ is SE-$\Rcal$ to $A$, then there is a nilpotent $N$
such that
$B$ is  SSE-$\Rcal$ to $A \oplus N$.
(See \cite{BoSc1} for further results, explanation and
context.)  A matrix is SE-$\Rcal$ to $(0)$
if and only if it is nilpotent. So, in particular,
 $N$ is SSE-$\Rcal$ to $(0)$ if and only if
the nilpotent matrix $N$ is trivial
as an element of $Nil_0(\Rcal )$.

\appendix

\section{Calculation of $\partial_{1}(\langle \epsilon, x+\epsilon \rangle)$}
This appendix contains the calculation of $\partial_{1}(\langle \epsilon, x+\epsilon \rangle)$. We let $I$ denote the ideal $(2)$ in $\mathbb{Z}[i]$, so that the map $\mathbb{Z}[i] \to \mathbb{F}_{2}[\epsilon]/(\epsilon^{2})$ given by $i \mapsto 1+\epsilon$ has kernel $I$. We also follow the notation in \cite[III.5.11]{WeibelBook}, so that, by definition,
$$\langle \epsilon, x + \epsilon \rangle = x_{ji}(-(x+\epsilon)(1-\epsilon x)^{-1})x_{ij}(-\epsilon)x_{ji}(x+\epsilon)x_{ij}((1-\epsilon x)^{-1}\epsilon)(h_{ij}(1-\epsilon x))^{-1}$$
where $i \ne j$, $x_{ij}, x_{ji}$ denote, as usual, generators in the Steinberg group, and for any unit $a$,
$$h_{ij}(a) = x_{ij}(a)x_{ji}(-a^{-1})x_{ij}(a)x_{ij}(-1)x_{ji}(1)x_{ij}(-1)$$
Steinberg relations give reduce $\langle \epsilon, x + \epsilon \rangle$ to
$$X := x_{ji}(-x-\epsilon-\epsilon x^{2})x_{ij}(-\epsilon)x_{ji}(x+\epsilon)x_{ij}(\epsilon)(h_{ij}(1-\epsilon x))^{-1}$$
Now $\partial_{1}$ is computed by composing up through the diagram, from bottom left to top right (see \cite[III.5.71]{WeibelBook} - we have introduced the $\phi_{i}$ names for the maps),
\[ \xymatrix{
& & GL(I) \ar^{\phi_{5}}[r] \ar^{\phi_{4}}[d] & K_{1}(\mathbb{Z}[i][x],(2)) \\
& St(\mathbb{Z}[i][x]) \ar^{\phi_{3}}[r] \ar^{\phi_{2}}[d] & GL(\mathbb{Z}[i][x]) & \\
K_{2}(\mathbb{F}_{2}[\epsilon,x]/(\epsilon^{2})) \ar[r] & St(\mathbb{F}_{2}[\epsilon,x]/(\epsilon^{2}) & & & }
\]
Now for simplicity we let $i=1, j=2$. We have $X \in St(\mathbb{F}_{2}[\epsilon,x]/(\epsilon^{2})$, and lifting $X$ up using $\phi_{2}$ and applying $\phi_{3}$ gives $\phi_{2}(YZ) = X$ where
$$Y := e_{21}(-x+1-i+(1-i)x^{2})e_{12}(1-i)e_{21}(x+i-1)e_{12}(i-1)$$
$$Z := e_{12}(1)e_{21}(-1)e_{12}(1)e_{12}((i-1)x-1)e_{21}(1+(i-1)x)e_{12}((i-1)x-1)$$
Lifting up via $\phi_{4}$ and applying $\phi_{5}$ gives the class of $YZ$ in $K_{1}(\mathbb{Z}[i][x],(2))$.
\bibliographystyle{plain}

\bibliography{mbssbib}

\end{document}